\documentclass{article}
\usepackage[utf8]{inputenc}

\usepackage{mathrsfs}
\usepackage{amsmath,amssymb}
\usepackage{amsthm,bm}
\usepackage[colorlinks,linkcolor=red,anchorcolor=blue,citecolor=green]{hyperref}
\usepackage{graphics,graphicx}
\usepackage{color}
\usepackage{enumerate}
\usepackage[title]{appendix}

\usepackage{subfigure}

\usepackage{caption}
\allowdisplaybreaks

\usepackage{tikz}
\usepackage{indentfirst}

\newtheorem{Thm}{Theorem}[section]
\newtheorem{Def}[Thm]{Definition}
\newtheorem{Lemma}[Thm]{Lemma}
\newtheorem{Coro}[Thm]{Corollary}

\newtheorem{Prop}[Thm]{Proposition}

\def\qed{\hfill \rule{4pt}{7pt}}

\def\pf{\noindent {\it{Proof.}\hskip 2pt}}

\title{On Negative Correlation of Arboreal Gas on Some Graphs
	\footnote{}
}
\author{
	Xiangyu Huang\footnote{YMSC, Tsinghua University, Beijing 100084, China. \emph{xiangyuhuang077@gmail.com}} 
}
\date{}

\usepackage[square,numbers]{natbib}

\begin{document}

\maketitle
	\begin{abstract}
		Arboreal Gas is a type of (unrooted) random forest on a graph, where the probability is determined by a parameter $\beta>0$ per edge. This model is essentially equivalent to acyclic Bernoulli bond percolation with a parameter $p=\beta/(1+\beta)$. Additionally, Arboreal Gas can be considered as the limit of the $q$-states random cluster model with $p=\beta q$ as $q\to 0$. A natural question arises regarding the existence and performance of the weak limit of Arboreal Gas as the graph size goes to infinity. The answer to this question relies on the negative correlation of Arboreal Gas, which is still an open problem. This paper primarily focuses on the negative correlation of Arboreal Gas and provides some results for specific graphs.\\
		
		%functionals of processes on graphs. Then we establish a novel dynamic programming equation associated with these functionals of the ORRW on finite graphs. Further, we give its solution by a variational representation. Based on this representation,
%we obtain the variational formula of the limit (as $n\to\infty$) of the functionals.
%By this variational formula, we verify the large deviation principle for empirical measures of the ORRW on finite graphs. In addition, we observe that the rate function is continuous and decreasing in $\delta$, and that it is not differentiable at $\delta=1$.

%Futhermore, for a class of stopping times, we derive the critical exponent for their exponential integrability precisely. Here, the stopping times contain the cover time and the hitting time. Moreover, we prove that the critical exponent is continuous, decreasing in $\delta$. Finally, we characterize the relationship between the limit (as $\delta\to0$) of the critical exponent and the structure of the graph.\\ %i.e. for $\delta$-ORRW there exists a $\alpha_c^1(\delta)$ such that $\E e^{\alpha C_E}<\infty$ for $\alpha<\alpha_c^1$ and $=\infty$ otherwise.
		%We believe this causes the difference between the rate function of ORRW and that of SRW.\\
		
	\end{abstract}

\section{Introduction to Arboreal Gas}
\noindent Let $G=(V,E)$ be a finite connected graph. If there exists an edge $e\in E$ connecting vertices $u,v\in V$, we say that $u,v$ are adjacent, denoted by $u\sim v$. We can also represent the edge $e$ as $uv$. A forest on $G$ is a subgraph that does not contain any cycle. Set $\mathcal{F}$ to be the collection of all forests on $G$. The \emph{Arboreal Gas} with parameter $\beta_e>0$ for each $e\in E$ is the measure on forests $F$  given by
\[
\mathbb{P}_\beta[F]:=\frac{1}{Z_\beta}\prod_{e\in F}\beta_e,\hskip 6mm Z_\beta:=\sum_{F\in\mathcal{F}}\prod_{e\in F}\beta_e.
\]
Specifically, if $\beta_e\equiv \beta$ for all edges $e$, the measure will be given by
\[
\mathbb{P}_\beta[F]:=\frac{1}{Z_\beta}\beta^{|F|},\hskip 6mm Z_\beta:=\sum_{F\in\mathcal{F}}\beta^{|F|},
\]
where $|F|$ stands for the number of edges in $F$.

Let $\mathbb{P}^{\rm perc}_{p}$ denote the probablity of Bernoulli bond percolation with parameter $p$. We can note that Arboreal Gas with a uniform parameter $\beta$ is equivalent to Bernoulli bond percolation with parameter $p_\beta:=\beta/(1+\beta)$ conditioned to be acyclic:
\[
\mathbb{P}^{\rm perc}_{p_\beta}[F|{\rm acyclic}]=\frac{p_\beta^{|F|}(1-p_\beta)^{|E|-|F|}}{\sum_{F\in \mathcal{F}}p_\beta^{|F|}(1-p_\beta)^{|E|-|F|}}=\frac{\beta^{|F|}}{\sum_{F\in\mathcal{F}}\beta^{|F|}}=\mathbb{P}_\beta[F].
\]
Another notable observation is that Arboreal Gas emerges as the limit of the $q$-states random cluster model as $q\to 0$, with $p=\beta q$, as mentioned in \cite{P2000}. The uniform forest model, also mentioned in \cite{P2000}, can be seen as one particular instance of Arboreal Gas, where $\beta=1$. It is important to note that this model is distinct from another known model referred to as \emph{uniform spanning forest} in the literature. Furthermore, another specific case is the uniform spanning tree, which arises as a limit of Arboreal Gas as $\beta\to\infty$.

An interesting phenomenon lies in the correlation between Arboreal Gas and hyperbolic spin systems. This kind of spin system is different from the classical spin systems with spherical symmetry, such as Ising model and classical Heisenberg model. Hyperbolic spin systems are also extensively studied in condensed matter physics due to their connection to the Anderson delocalisation-localisation transition of random Sch\"{o}dinger operators and related matrix models, as discussed in \cite{E1983,SW1980,W1979}. Researchers are interested in the existence of phase transitions in the hyperbolic spin systems. Considering Arboreal Gas, Bauerschmidt, Crawford, Helmuth and Swan presented a magic formula in \cite{BCHS2021} that elegantly describe its probability by the intergral in the hyperbolic spin system $\mathbb{H}^{0|2}$, which may also be seen in \cite{BH2021,CJSSS2004}.

The subcritical percolation and percolating phase transitions in Arboreal Gas have garnered significant attention.
In the seminal works \cite{BCS2009,LP1992,MY2017}, it was established that Arboreal Gas, with a parameter $\beta=\alpha/N$ for fixed $\alpha$, undergoes a phase transition on the complete graph $K_N$ with $N$ vertices. This phase transition was also demonstrated by Bauerschmidt, Crawford, Helmuth, and Swan in their recent study \cite{BCHS2021}. Additionally, they proved the polynomial decay of connectivity probability from the origin to other vertices in all subgraphs of the lattice $\mathbb{Z}^2$. Based on this polynomial decay, and assuming for the existence of a weak limit of probabilities on asymptotic graphs of $\mathbb{Z}^2$, they demonstrated the absence of an infinite tree almost surely. On torus $\Lambda_N=\mathbb{Z}^d/ L^N\mathbb{Z}^d$ with $L$ fixed and $N\to\infty$, Bauerschmidt, Crawford and Helmuth gave an asymptotic estimate of connectivity probability from the origin to other vertices for $d\ge 3$ in \cite{BCH2021}. Recently, Halberstam and Hutchcroft verified the uniqueness of infinite tree in $\mathbb{Z}^d$ for $d=3,4$ for any translation-invariant Arboreal Gas Gibbs measures in \cite{HH2023}.

For any edge set $S\subset E$, let $\mathbb{P}_\beta[S]$ be the probability that all edges in $S$ belong to the forest. We simply write $\mathbb{P}_\beta[S]$ as $\mathbb{P}_\beta[e_1e_2\dots e_n]$ for $S=\{e_1,e_2,\dots,e_n\}$. For two vertices $u,v\in V$, denote by $\mathbb{P}_\beta[u\leftrightarrow v]$ the probability that $u,v$ are in a same tree of the forest. The \emph{negative correlation} of Arboreal Gas should be expressed as
\begin{equation}
\mathbb{P}_\beta[e_1e_2]\le \mathbb{P}_\beta[e_1]\mathbb{P}_\beta[e_2],\hskip 6mm {\rm for\ any\ }e_1,e_2\in E.\label{NC}
\end{equation}

Until now, this question is still open. A weaker inequality has been proved recently in \cite{BH2022} by Br\"{a}nd{\'e}n and Huh. They showed $\mathbb{P}_\beta[e_1e_2]\le 2\mathbb{P}_\beta[e_1]\mathbb{P}_\beta[e_2]$ by using the Lorentzian signature. Once (\ref{NC}) holds, the existence of the weak limit of Arboreal Gas Gibbs measures on increasing asymptotic graphs $G_n\uparrow \mathbb{Z}^2$ as $n\to\infty$ was answered in \cite{BCHS2021}. Under this weak limit on $\mathbb{Z}^2$, all trees should be finite almost surely, see Corollary 1.4 in \cite{BCHS2021}.

In this paper, we focus on the negative correlation of Arboreal Gas on finite connected graph $G$. One way to consider this question is to realize the relevance between Arboreal Gas and the hyperbolic spin system $\mathbb{H}^{0|2}$. The negative correlation of Arboreal Gas can be implied by the monotonicity of a hyper-integral on $\mathbb{H}^{0|2}$, see section 5.2 in \cite{BH2021}. However, the monotonicity of hyper-integral is hard to be verified. In our paper, we consider Arboreal Gas directly. Here is the outline of this paper. In Section \ref{sec main}, we show our main results, and we give some preliminaries in Section \ref{sec pre}. In Section \ref{sec nc adjacent edges}, we show the negative correlation for adjacent edges and sufficiently large $\beta$. In Section \ref{sec nc complete graph}, we show the negative correlation of complete graphs $K_n$ with  $n$ vertices for sufficiently large $n$ and $\beta$. In Section \ref{sec nc simplify graph}, we give the proof of the equivalence of the negative correlation on any graph and on its simplified version.

\section{Negative Correlations and Main Results}\label{sec main}
\noindent   Here are our main theorems.

\begin{Thm}\label{thm nc adjacent edges}
	Consider Arboreal Gas with parameter $\beta$ on finite connected graph $G$. For sufficiently large $\beta$ and any adjacent distinct edges $e_1,e_2\in E$,
	\[
	\mathbb{P}_\beta[e_1e_2]\le \mathbb{P}_\beta[e_1]\mathbb{P}_\beta[e_2].
	\]
\end{Thm}

\begin{Thm}\label{thm nc complete graph}
	Consider Arboreal Gas with parameter $\beta$ on complete graph $K_n$ with $n$ vertices, where $n$ is large enough. For sufficiently large or sufficiently small $\beta$ and any distinct edges $e_1,e_2\in E$,
	\[
	\mathbb{P}_\beta[e_1e_2]\le \mathbb{P}_\beta[e_1]\mathbb{P}_\beta[e_2].
	\]
\end{Thm}

We verify those two theorems by regarding the probability $\mathbb{P}_\beta[e_1e_2], \mathbb{P}_\beta[e_1], \mathbb{P}_\beta[e_2]$ as polynomials of $\beta$. A fact is, the negative correlation holds for uniform spanning trees on finite connected graphs, see  section 4.2 in \cite{LP2017}. An important observation is that the term with highest degree in the polynomials of $\beta$ mentioned before corresponds to the uniform spanning trees. This is our basic idea to show the negative correlation for sufficiently large $\beta$.\\ 

Besides, we have another theorem to simplify the structure of graphs when we consider the negative correlation. Before we give this theorem, we first give some definitions. 

%\begin{Thm}\label{torus G, small beta, adjacent edge}
%Consider Arboreal Gas with parameter $\beta$ on some torus on $\mathbb{Z}^2$. For $\beta\le1/16$ and any distinct edges $e_1,e_2\in E$ in a same square with length $1$,
%\[
%\mathbb{P}_\beta[e_1e_2]\le \mathbb{P}_\beta[e_1]\mathbb{P}_\beta[e_2].
%\]
%\end{Thm}

%This theorem also comes from the analyse of the probability $\mathbb{P}_\beta[e_1e_2], \mathbb{P}_\beta[e_1], \mathbb{P}_\beta[e_2]$. The difference is that we turn it into the analyse of the structure of the graph.

\begin{Def}
	For $e\in E$, we call it pivotal for $G$ if extreme points of $e$ are not connected in $G\setminus e:=(V,E\setminus\{e\})$.
\end{Def}

\begin{Def}\label{def delete 2-d vertices}
	For a graph $G$ with no pivotal edges, we give graph $\tilde{G}=(\tilde{V},\tilde{E})$ by induction:
	\begin{enumerate}
		\item Set $G_0=G$.
		\item If there exists some vertex $u$ in $V_n$ having degree $2$, we assume $v_1,v_2$ are adjacent to $u$. Then we set $V_{n+1}=V_n\setminus \{u\}$. We give $E_{n+1}$  by adding a new edge $e_0$ connecting $v_1,v_2$ on $E_n\setminus\{uv_1,uv_2\}$. Define $f_{n+1}: E_n\to E_{n+1}$ such that $f_{n+1}(uv_i)=e_0$ for $i=1,2$ and $f_{n+1}(e)=e$ for $e\neq uv_1,uv_2$.
		\item Stop if there is no vertex with degree $2$.
	\end{enumerate}
	Further, we give a map $f:E\to \tilde{E}$ by $f=f_n\circ f_{n-1}\circ\dots\circ f_1$, where $n$ is the number of steps to get $\tilde{E}$.
\end{Def}

%\begin{Def}
%For two edges $e_1,e_2\in E$, they are in a same equivalent class if there is a path connecting $e_1,e_2$ such that all vertices in this path have degree $2$.
%\end{Def}

\begin{Def}\label{def delete multiple edges}
	For a graph $G$, we give a simple graph $G'=(V',E')$ by setting $V'=V$ and letting $u,v\in V'$ adjacent in $G'$ if they are adjacent in $G$. Further, we give a map $g:E\to E'$ with $g(e)=uv$ if the end points of $e$ are $u,v$.
\end{Def}

Here we give our main theorem.

\begin{Thm}\label{NG for same edges}
	Consider Arboreal Gas with parameter $\beta_e$ for each $e\in E$ on finite connected graph $G$. Give $\hat{G}$ by deleting all pivotal edges in $G$, then the negative correlation on $G$ is equivalent to that on each component of $g\circ f(\hat{G})$. 
\end{Thm}

\begin{Coro}
	Consider Arboreal Gas with parameter $\beta_e$ for each $e\in E$ on finite connected graph $G$. For any distinct edges $e_1,e_2\in E$, if $e_1$ or $e_2$ is pivotal, or $g\circ f(e_1)=g\circ f(e_2)$, then we have
	\[
	\mathbb{P}_\beta[e_1e_2]\le \mathbb{P}_\beta[e_1]\mathbb{P}_\beta[e_2].
	\]
\end{Coro}

This corollary immediately comes out by Theorem \ref{NG for same edges}. We can also obtain the negative correlation on some specific graphs.

\begin{Coro}
	Consider Arboreal Gas with parameter $\beta_e$ for each $e\in E$. . For any integer $d\ge0$, the negative correlation holds on ladder $\mathbb{L}_d:=\{1,2,\dots,d \}\times \{0,1\}$.
\end{Coro}

\pf Note the fact that $g\circ f(\mathbb{L}_d)=\mathbb{L}_{d-2}$ for $d\ge2$, where $\mathbb{L}_0$ is a unique vertex, and $\mathbb{L}_1$ is a line segment. By induction, the negative correlation on $\mathbb{L}_d$ can be implied by that on $\mathbb{L}_0$ and $\mathbb{L}_1$, which completes the proof.\qed

\section{Preliminaries}\label{sec pre}

\noindent First we show an equivalent condition of negative correlation.

\begin{Def}
	For disjoint $S_1,S_2\subset E$, define
	\[
	\mathbb{P}_\beta[S_1\bar{S_2}]=\mathbb{P}_\beta[S_1\subset F, S_2\cap F=\emptyset ].
	\]
	Given a measure $\mu$ defined on $G$ by
	\[
	\mu[S_1\bar{S_2}]=\mathbb{P}_\beta[S_1\bar{S_2}]\cdot Z_\beta / \prod_{e\in S_1}\beta_{e}.
	\]
	Specifically, set $\mu[1]=Z_\beta$.
\end{Def}

\begin{Prop}
	The following two conditions are equivalent:
	\begin{enumerate}
		\item $\mathbb{P}_\beta[e_1e_2]\le \mathbb{P}_\beta[e_1]\mathbb{P}_\beta[e_2]$ for all distinct $e_1,e_2\in E$;
		\item $\mathbb{P}_\beta[S_1S_2]\le \mathbb{P}_\beta[S_1]\mathbb{P}_\beta[S_2]$ for all disjoint $S_1,S_2\subset E$.
	\end{enumerate}
\end{Prop}

\pf

\noindent $1\Leftarrow 2$: This result immediately comes out by setting $S_i=e_i$ for $i=1,2$.

\noindent $1\Rightarrow 2$: First we prove the equivalence of condition $1$ and decreasing of $\mathbb{P}_\beta[e_0]$ in each $\beta_e$ for all $e_0\in E$ and $e\neq e_0$. Note that
\begin{align}
1 &\Leftrightarrow  \mathbb{P}_\beta[e_1e_2]\mathbb{P}_\beta[\bar{e_1}\bar{e_2}]\le \mathbb{P}_\beta[e_1\bar{e_2}]\mathbb{P}_\beta[\bar{e_1}e_2]\nonumber\\
&\Leftrightarrow  \mu[e_1e_2]\mu[\bar{e_1}\bar{e_2}]\le\mu[e_1\bar{e_2}]\mu[\bar{e_1}e_2]\nonumber\\
&\Leftrightarrow  \frac{\beta_{e_1}\mu[e_1\bar{e_2}]}{\beta_{e_1}\mu[e_1\bar{e_2}]+\mu[\bar{e_1}\bar{e_2}]}\le \frac{\beta_{e_1}\mu[e_1e_2]}{\beta_{e_1}\mu[e_1e_2]+\mu[\bar{e_1}e_2]}\nonumber\\
&\Leftrightarrow  \frac{\beta_{e_2}\beta_{e_1}\mu[e_1e_2]+\beta_{e_1}\mu[e_1\bar{e_2}]}{\beta_{e_2}(\beta_{e_1}\mu[e_1e_2]+\mu[\bar{e_1}e_2])+\beta_{e_1}\mu[e_1\bar{e_2}]+\mu[\bar{e_1}\bar{e_2}]} {\rm\ is\ decreasing\ in\ } \beta_{e_2}.\label{1 in ng equiv}
\end{align}
By $\beta_{e_2}\mu[\bar{e_1}e_2]+\mu[\bar{e_1}\bar{e_2}]=\mu[\bar{e_1}]$ and $\beta_{e_2}\mu[e_1e_2]+\mu[e_1\bar{e_2}]=\mu[e_1]$,
\[
(\ref{1 in ng equiv})=\mathbb{P}_\beta[e_1].
\]
By the arbitrary of $e_1,e_2$, we obtain that (\ref{1 in ng equiv}) is equivalent to the decreasing of $\mathbb{P}_\beta[e_0]$ in each $\beta_e$ for all $e_0\in E$ and $e\neq e_0$.\\

Secondly, we prove that the decreasing of $\mathbb{P}_\beta[e_0]$ in each $\beta_e$ for all $e_0\in E$ and $e\neq e_0$ deduces condition $2$, which may complete the proof of $1\Rightarrow 2$. To show this property, we claim that

{\begin{center}$\mathbb{P}_\beta[\cdot|S]$ equals the limit of $\mathbb{P}_\beta[\cdot]$ as $\beta_e\to\infty$ for all $e\in S\subset E$.\end{center}}

By the decreasing of $\mathbb{P}_\beta[e_0]$ in each component of $\beta$ except $\beta_{e_0}$ for all $e_0\in E$, we know that for $S\subset E$,
\[
\mathbb{P}_\beta[e|S]\le\mathbb{P}_\beta[e] {\rm\ for\ all\ edge\ }e\notin S.
\]
For disjoint $S_1,S_2\subset E$, Applying this to $\mathbb{P}[\cdot|S_1]$, we obtain 
\[
\mathbb{P}_\beta[e|S_1S_2]\le\mathbb{P}_\beta[e|S_1] {\rm\ for\ all\ edge\ }e\notin S_1,S_2.
\]
This inequality is equivalent to
\[
\mathbb{P}_\beta[S_2|eS_1]\le \mathbb{P}_\beta[S_2|S_1].
\]
BY induction of $S_1$, we verify that
\[
\mathbb{P}[S_2|S_1]\le \mathbb{P}[S_2],
\]
which implies condition $2$.

Finally we prove our claim by setting
\[
\beta_e^{x,S}=\left\{\begin{aligned}
&\beta_e, &e\notin S,\\
&x\cdot \beta_e, &e\in S.
\end{aligned}
\right.
\]
Then $\mathbb{P}_{\beta^{x,S}}[F]=\prod_{e\in F}\beta_e^{x,S}/\sum_{F'\in \mathcal{F}}\prod_{e\in F'}\beta_e^{x,S}$, which converges to (as $x\to\infty$)
\[
\frac{\prod_{e\in F\setminus S}\beta_e \textbf{1}_{\{F\supset S\}}}{\sum_{F'\supset S}\prod_{e\in F'\setminus S}\beta_e \textbf{1}_{\{F'\supset S\}}}.
\]
By $\mathbb{P}_\beta[F|S]\propto \prod_{e\in F\setminus S}\beta_e \textbf{1}_{\{F\supset S\}}$,
we show that $\mathbb{P}_{\beta^{x,S}}[\cdot]$ converges to $\mathbb{P}_\beta[\cdot|S]$ in weak topology. Further, by the decreasing of $\mathbb{P}_{\beta^{x,S}}[\cdot]$ in $x$, we conclude that $\mathbb{P}_\beta[\cdot|S]\le \mathbb{P}_\beta[\cdot]$.   \qed\\

%%%%electrical network

Next we show the negative correlation of uniform spanning trees. Here we give the definition of the uniform spanning trees.

\begin{Def}
A spanning tree $T$ on $G$ is a subgraph of $G$ containing all vertices in $V$ and being a tree. The \emph{uniform spanning tree} is a uniform measure on all spanning trees on $G$.
\end{Def}

In our next part, denote by $\mathbb{P}_{\infty}^{G}$ the probability of uniform spanning trees on graph $G$. In fact, uniform spanning trees can be generated by simple random walks. This method is called \emph{Wilson's algorithm}, see the details of which in Theorem 4.1 of \cite{LP2017}. To show the negative correlation of uniform spanning trees, we use this method and another tool called \emph{electrical network}. This tool shows a nice relationship between the probability theory and potential theory. Here we define the electrical network.\\

For a finite connected graph $G$ with conductance $c_e$ on each edge, we build an electrical network. Given  two disjoint vertex sets $S$ and $T$ to be source and sink. Define potential difference $\phi(\vec{e})$ and current $i(\vec{e})$ on each direct edge $\vec{e}$.
For adjacent vertices $u,v$,
\[
\phi(\overrightarrow{uv})=-\phi(\overrightarrow{vu}),\hskip 3mm i(\overrightarrow{uv})=-i(\overrightarrow{vu}).
\]
Here we normally take
\[
\sum_{u\sim v}i(\overrightarrow{uv})\left\{
\begin{aligned}
&>0, &v\in S,\\
&<0, &v\in T.
\end{aligned}
\right.
\]
Further, the potential difference and current satisfies the following three laws, see \cite{G2011}.\\

\noindent \textbf{Kirchhoff's potential law.} For any cycle $v_1v_2\cdots v_nv_{n+1}$ with $v_{n+1}=v_1$,
\[
\sum_{i=1}^{n}\phi(\overrightarrow{v_iv_{i+1}})=0.
\]
\noindent \textbf{Kirchhoff's current law.} For any $v\notin S\cup T$,
\[
\sum_{u\sim v}i(\overrightarrow{uv})=0.
\]
\noindent \textbf{Ohm's law.} For any edge $e=uv$,
\[
i(\vec{e})c(e)=\phi(\vec{e}).
\]
Kirchhoff's potential law is equivalent to the existence of the function $\phi$ on vertex set such that $\phi(\overrightarrow{uv})=\phi(v)-\phi(u)$ for adjacent vertices $u,v$. Here $\phi$ is also called potential function. Next we define the energy of current $i$ by
\[
E(i)=\frac{1}{2}\sum_{u,v\in V}i^2(\overrightarrow{uv})/c(uv).
\]
Define the unit currenct flow $i$ to be such a flow that $\sum_{v\in S,u\sim v}i(\overrightarrow{uv})=1$. Denote by 
\[
R_{\rm eff}(G,c)=E(i)\ {\rm for\ a\ unit\ current\ }i
\]
the effective resistance of the electrical network $G$ with conductance $c$.
Here $R_{\rm eff}(G,c)$ is equal to $\phi(v)$ for $v\in S$ if we set $\phi(u)=0$ for $u\in T$, where the current is unit.\\

For some finite connected graph $G$, we have the following equation for its spanning trees.

\begin{Prop} [Kirchhoff's Effective Resistance Formula, \cite{LP2017}]\label{Prop R-eff}
	Let $T$ be the spanning tree of a finite connected graph $G$, and $e$ be some edge of $G$ with endpoints $e^-,e^+$, then
	\[
	\mathbb{P}_{\infty}^G[e\in T]=\mathbb{P}_{e^-}[1st\ hits\ e^+\ by\ travelling\ along\ e]=i(\overrightarrow{e^-e^+})=c(e)R_{\rm eff}(G),
	\]
	where $i$ is unit current flow from $e^-$ to $e^+$.
\end{Prop}

By Kirchhoff's Effective Resistance Formula, we deduce the negative correlation of uniform spanning trees. Before that, we give a lemma to show that the effective resistance decreases if the resistance on some edge decreases by the following lemma.

\begin{Lemma}{\rm (Rayleigh principle, \cite{G2011})}\label{lemma R decreases in r}
	Consider an electrical network with unit current flow $i$ from $s$ to $t$ and conductance $c_e$ on each edge $e\in E$. For any edge $e_0\in E$, if $i(\vec{e_0})\neq 0$, then the effective resistance $R_{\rm eff}(c)$ strictly decreases in $c(e_0)$. Otherwise, $R_{\rm eff}(c)$ is a constant in $c(e_0)$.
\end{Lemma}

\noindent \emph{Proof of Lemma \ref{lemma R decreases in r}.} For the unit current flow $i$ and conductance $c$, we consider the energy
\[
E_c(i)=\frac{1}{2}\sum_{u,v\in V}i^2(\overrightarrow{uv})/c(uv),
\]
which is equal to $R_{\rm eff}(c)$.
Set $c'$ to be another conductance with
\[
c'(e)\left\{\begin{aligned}
&>c(e_0), &e=e_0,\\
&=c(e), &e\neq e_0.
\end{aligned}
\right.
\]
If $i(\vec{e_0})\neq 0$, we may observe that $E_{c'}(i)<E_c(i)$. Since the unit current flow $i'$ on the electrical network with conductance $c'$ has a smaller energy than flow $i$, we obtain
\[
R_{\rm eff}(c')=E_{c'}(i')\le E_{c'}(i)<E_c(i).
\]

If $i(\vec{e_0})=0$, we have $E_{c'}(i)=E_{c}(i)$. Here we prove that $i$ is the unit current flow on the electrical network with conductance $c'$. Consider the potential $\phi$ on the electrical network with conductance $c$, then for any edge $e\in E$,
\[
i(\vec{e})c(e)=\phi(\vec{e}).
\]
This implies that $\phi(\vec{e_0})=0$. On the electrical network with conductance $c'$, $\phi$ and $i$ satisfy Kirchhoff's potential law and Kirchhoff's current law. As to Ohm's law, for $e\neq e_0$,
\[
i(\vec{e})c'(e)=i(\vec{e})c(e)=\phi(\vec{e}).
\]
For $e=e_0$,
\[
i(\vec{e_0})c'(e_0)=0=\phi(\vec{e}).
\]
Therefore, $i$ is a unit current flow on the electrical network with conductance $c'$, which implies
\[
R_{\rm eff}(c')=E_{c'}(i)=E_{c}(i).
\]
\qed

\begin{Prop}{\rm \cite{LP2017}}\label{prop NC of UST}
	The negative correlation holds for uniform spanning trees on finite connected graph $G$.
\end{Prop}

\pf We only need to verify that for any two distinct edges $e_1,e_2$,
\[
\mathbb{P}_{\infty}^G[e_1e_2]\le \mathbb{P}_{\infty}^G[e_1]\mathbb{P}_{\infty}^G[e_2].
\]
This inequality is equivalent to
\[
\mathbb{P}_{\infty}^{G/e_2}[e_1]=\frac{\mathbb{P}_{\infty}^G[e_1e_2]}{\mathbb{P}_{\infty}^G[e_2]}\le \mathbb{P}_{\infty}^G[e_1],
\]
where $G/e_2$ stands for the contraction of $G$ by removing $e_2$ and identifying its endpoints. By Proposition \ref{Prop R-eff}, this inequality is equivalent to
\begin{equation}
R_{\rm eff}(G/e_2)\le R_{\rm eff}(G),\label{eq resistance of trees}
\end{equation}
where the unit current flow is from $e_1^-$ to $e_1^+$, and $e_1^-,e_1^+$ are endpoints of $e_1$. Actually, \eqref{eq resistance of trees} holds since $G/e_2$ can be regarded as graph $G$ with resistance $0$ on $e_2$, leading to a smaller effective resistance.\qed

\section{Proof of Theorem \ref{thm nc adjacent edges}}\label{sec nc adjacent edges}

\noindent For finite connected graph $G$, set $T[1]$ to be the number of spanning trees on $G$. Then for each spanning tree, it has the probability $1/T[1]$. For edge sets $S_1$ and $S_2$, let $T[S_1\bar{S_2}]$ be the number of all spanning trees containing $S_1$ and disjoint to $S_2$. Consider the Arboreal Gas with parameter $\beta$. Observe that the highest term of $\mu[S_1\bar{S_2}]$ in $\beta$ for some edge $e\in F$ should be
\[
T[S_1\bar{S_2}]\beta^{|V|-1}.
\] 
If we consider the negative correlation of Arboreal Gas, i.e.,
\[
\mu[e_1]\mu[e_2]-\mu[e_1e_2]\mu[1]\ge0,
\]
we pay attention to the highest term of the left hand side of this inequality if $\beta$ is large enough.
Here the highest term is
\[
(T[e_1]T[e_2]-T[e_1e_2]T[1])\beta^{2(|V|-1)}.
\]
By Proposition \ref{prop NC of UST}, this value should be non-negative. What we want is to show
\[
T[e_1]T[e_2]-T[e_1e_2]T[1]>0,
\]
by which we conclude the negative correlation for $\beta$ large enough. However, this does not always hold, such as that on trees. To prove Theorem \ref{thm nc adjacent edges}, for adjacent edges $e_1$ and $e_2$, we consider in two situations:
\begin{itemize}
	\item[(a)] there is a simple cycle crossing both $e_1$ and $e_2$;
	\item[(b)] there is no simple cycle crossing both $e_1$ and $e_2$.
\end{itemize}

In Case (a), if we start a unit current flow $i$ from $e_1^-$ to $e_1^+$, then $i(\vec{e_2})$ is not vanishing.

\begin{Lemma}\label{lemma R decreases in loop}
	For two adjacent edges $e_1,e_2\in E$, set $i$ to be the unit current flow from $e_1^-$ to $e_1^+$. If there is a simple cycle on $G$ crossing both $e_1$ and $e_2$, then $i(\vec{e_2})$ is not vanishing.
\end{Lemma}

In Case (b), we prove that $e_1,e_2$ can be separated into two parts of the graph $G$, where the intersection of these two parts has only one vertex. By this property, the probability of Arboreal Gas on $G$ can be expressed as the product of the probabilities of Arboreal Gas on these two parts.

\begin{Lemma}\label{lemma separate graph into 2 parts}
	For two adjacent edges $e_1,e_2\in E$, assume that $e_i^-,e_i^+$ are two endpoints of edge $e_i$ for $i=1,2$, where $e_1^+=e_2^+$. If there is no simple cycle on $G=(V,E)$ crossing both $e_1$ and $e_2$, then there exist subgraphs $G_i=(V_i,E_i)$ for $i=1,2$ such that
	\begin{itemize}
		\item $e_i\in E_i$ for $i=1,2$;
		\item $V_1\cap V_2=\{e_1^+\}$, $E_1\cap E_2=\emptyset$;
		\item $V_1\cup V_2=V$, $E_1\cup E_2=E$.
	\end{itemize}
\end{Lemma}

Now we prove Theorem \ref{thm nc adjacent edges}.\\

\noindent \emph{Proof of Theorem \ref{thm nc adjacent edges}.} 

\noindent \textbf{(a).} If there is a simple cycle crossing both $e_1$ and $e_2$, by Lemmas \ref{lemma R decreases in loop} and \ref{lemma R decreases in r}, we obtain
\[
R_{\rm eff}(G/e_2)<R_{\rm eff}(G)
\]
for the unit current flow from $e_1^-$ to $e_1^+$.
Then the highest term in $\mu[e_1]\mu[e_2]-\mu[e_1e_2]\mu[1]$ is positive, i.e.,
\[
 (T[e_1]T[e_2]-T[e_1e_2]T[1]) \cdot\beta^{2|V|-2}=(R_{\rm eff}(G)-R_{\rm eff}(G/e_2))\cdot T[e_2]T[1]\beta^{2|V|-2}>0,
\]
where the equality comes from Proposition \ref{Prop R-eff}.
This implies the negative correlation for $e_1,e_2$, i.e., $\mu[e_1]\mu[e_2]-\mu[e_1e_2]\mu[1]\ge0$, for sufficiently large $\beta$.\\

\noindent \textbf{(b).} If there is no simple cycle crossing both $e_1$ and $e_2$ (assume $e_1^+=e_2^+$), by Lemma \ref{lemma separate graph into 2 parts}, we separate $G$ into two parts $G_1=(V_1,E_1)$ and $G_2=(V_2,E_2)$ with
\begin{itemize}
	\item $e_i\in E_i$ for $i=1,2$;
	\item $V_1\cap V_2=\{e_1^+\}$, $E_1\cap E_2=\emptyset$;
	\item $V_1\cup V_2=V$, $E_1\cup E_2=E$.
\end{itemize}
We claim that
\begin{equation*}
\text{there is no simple cycle crossing both edges in } E_1 \text{ and } E_2.
\end{equation*}
By this claim, we obtain that the subgraph $F$ in $G$ is a forest is equivalent to that $F|_{G_i}$ are forests for $i=1,2$, where $F|_{G_i}=(V(F)\cap V_i, E(F)\cap E_i)$. Therefore, for $i=1,2$, if we denote by $\mu_i$ the measure for Arboreal Gas on $G_i$, i.e., the probability multiplying its partition function, then
\begin{align*}
\mu[e_1\bar{e_2}]=\mu_1[e_1]\mu_2[\bar{e_2}],\\
\mu[\bar{e_1}e_2]=\mu_1[\bar{e_1}]\mu_2[e_2],\\
\mu[e_1e_2]=\mu_1[e_1]\mu_2[e_2],\\
\mu[\bar{e_1}\bar{e_2}]=\mu_1[\bar{e_1}]\mu_2[\bar{e_2}].
\end{align*}
By these properties, we deduce that
\[
\mu[e_1\bar{e_2}]\mu[\bar{e_1}e_2]-\mu[e_1e_2]\mu[\bar{e_1}\bar{e_2}]=0,
\]
which implies the negative correlation for $e_1,e_2$.\\

Now we prove the claim by contradiction. Otherwise, there is a simple cycle $C=u_1u_2\dots u_nu_{n+1}$ with $u_{n+1}:=u_1$ such that $u_{k_i}u_{k_i+1}\in E_i$ for $i=1,2$ and some $k_1\neq k_2$. Without loss of generality, assume $1\le k_1<k_2\le n$. Since $E_1\cup E_2=E$, there should be two distinct integers $n_1,n_2$ with $k_1\le n_1<k_2$ and $n_2<k_1$ or $n_2\ge k_2$ such that $u_{n_1}u_{n_1+1}, u_{n_2+1}u_{n_2+2}\in E_1$ and $u_{n_1+1}u_{n_2+2}, u_{n_2}u_{n_2+1}\in E_2$. This implies $u_{n_1+1},u_{n_2+1}\in V_1\cap V_2$, i.e., $u_{n_1+1}=u_{n_2+1}=e_1^+$, which is contradict to $n_1\neq n_2$ and that $C$ is a simple cycle.
\qed\\

Finally we prove Lemmas \ref{lemma R decreases in loop} and \ref{lemma separate graph into 2 parts}.\\

\noindent \emph{Proof of Lemma \ref{lemma R decreases in loop}.} Without loss of generality, assume $e_1^-$ is a endpoint of $e_2$.
We prove by contradiction. If $i(\vec{e_2})=0$, we consider the simple cycle $\pi$ crossing both $e_1$ and $e_2$, and set
\[\pi=e_1^- u_1u_2\dots u_n e_1^+ e_1^-.\]
Note that the potentials $\phi$ on $e_1^-,e_1^+$ are different. Precisely, $\phi(e_1^-)>\phi(e_1^+)=0$, which implies $i(\overrightarrow{e_1^- e_1^+})> 0$. By Kirchhoff's potential law and Ohm's law, we know that there should be some edge $e\neq e_1$ with $i(\vec{e})\neq 0$. Let $m:=\inf\{ k\ge 1: i(\overrightarrow{u_ku_{k+1}})\neq 0 \}$, where $u_{n+1}:=e_1^+$. By Kirchhoff's current law, we deduce that there is a vertex $v_1$ adjacent to $v_0:=u_m$ such that $i(\overrightarrow{v_0v_1})<0$. By Kirchhoff's potential law, $v_1\neq e_1^+$. Otherwise, potentials on the cycle $v_0v_1e_1^-u_1\dots u_{m-1}v_0$ do not obey Kirchhoff's potential law.
Now we construct a chain by induction:

For integer $k\ge1$, assume there is a chain $v_0v_1\dots v_k$ with $i(\overrightarrow{v_{j-1}v_{j}})<0$ for $j=1,\dots,k$ and $v_{j_1}\neq v_{j_2}\neq e_1^+$ for distinct $j_1,j_2=1,\dots,k$. By Kirchhoff's current law, we can find a vertex $v_{k+1}$ adjacent to $v_k$ with $i(\overrightarrow{v_{j-1}v_{j}})<0$. Furthermore, $v_{k+1}\neq v_j,e_1^+$ for $j=0,\dots,k$. Otherwise, if $v_{k+1}=e_1^+$, potentials on the cycle $v_0v_1\dots  v_ke_1^-u_1\dots u_{m-1}v_0$ do not obey Kirchhoff's potential law. if $v_{k+1}=v_j$ for some $j=0,1,\dots,k$, potentials on the cycle $v_jv_{j+1}\dots v_kv_j$ do not obey Kirchhoff's potential law.

Thus we construct a chain $v_0\dots v_{|V|}$ with $v_{j_1}\neq v_{j_2}$ for distinct $j_1,j_2=0,\dots,|V|$. That is, we find $|V|+1$ different vertices, which is contradict to that the number of all vertices of $G$ is $|V|$.
\qed\\

\noindent \emph{Proof of Lemma \ref{lemma separate graph into 2 parts}.} 
Without loss of generality, assume $G$ is connected. For $i=1,2$, set $G_i=(V_i, E_i)$ to be the graph induced by the vertex set
\[V_i:=\{e_i^+,e_i^-\}\cup\{ u\in V: {\rm there\ is\ a\ simple\ path\ from\ }e_i^-\ {\rm to\ }u\ {\rm not\ visiting\ }e_i^+ \}.\]
Set $G_3=(V_3,E_3)$ to be the graph induced by
\[
V_3:=\{e_1^+\}\cup [V\setminus(V_1\cup V_2)].
\]
Now we prove that $G_1,G_2\cup G_3$ satisfy conditions in Lemma \ref{lemma separate graph into 2 parts} by showing the following properties.
\begin{itemize}
	\item \underline{$e_i\in E_i$ for $i=1,2$.}
	
	Since $e_i^-,e_i^+\in V_i$ for $i=1,2$, $e_i\in E_i$.
	
	\item \underline{$V_i\cap V_j=\{e_1^+\}$ and $E_i\cap E_j=\emptyset$ for distinct $i,j=1,2,3$.}
	
	If there is an edge $e\in E_i\cap E_j$ for some distinct $i,j=1,2,3$, then two endpoints of $e$ belong to $V_i$ and $V_j$, contradict to $V_i\cap V_j=\{e_1^+\}$. Thus we only need to show $V_i\cap V_j=\{e_1^+\}$ for distinct $i,j=1,2,3$.
	
	For any $i=1,2$, $V_i\cap V_3=\{e_1^+ \}$ by the definition of $V_3$. We then show $V_1\cap V_2=\{e_1^+\}$.
	
	For any vertex $v\in V_1$ with $v\neq e_1^+$, set $\pi=e_1^-u_1\dots u_n v$ to be the simple path from $e_1^-$ to $v$ not visiting $e_1^+$. Then $e_2^-\notin \pi$, otherwise there is some $k$ with $u_k=e_2^-$ and there will be a simple cycle $e_1^+u_1\dots u_k e_1^+$ crossing $e_1,e_2$, contradict to the condition in Lemma \ref{lemma separate graph into 2 parts}.
	
	We prove $v\notin V_2$ by contradiction. If $v\in V_2$, $v\in\pi$ implies $v\neq e_2^-$. Therefore, there is a simple path $\pi'=v_0 v_1\dots v_{n'} e_2^-$ from $v$ to $e_2^-$ not visiting $e_1^+$, where $v_0:=v$. Set $k':=\sup\{ k\le n':v_k\in \pi \}$. Assume $v_{k'}=u_{k}$ for some $k$, where $u_0:=e_1^-$ and $u_{n+1}:=v$. Then $e_1^+ u_0\dots u_{k}v_{k'+1}\dots v_{n'} e_2^- e_1^+$ is a simple cycle crossing $e_1,e_2$, contradict to the condition in Lemma \ref{lemma separate graph into 2 parts}.
	
	Since we choose $v\in V_1$ arbitrarily, we obtain $V_1 \cap V_2=\{e_1^+\}$.
	
	\item \underline{$V_1\cup V_2\cup V_3=V$, $E_1\cup E_2\cup E_3=E$.}
	
	By definition of $V_3$, we immediately obtain $V_1\cup V_2\cup V_3=V$. Since $V_i\cap V_j=\{e_1^+\}$ for distinct $i,j=1,2,3$, we show $E_1\cup E_2\cup E_3=E$ by proving that there is no edge connecting $V_i,V_j$ except those edges with endpoint $e_1^+$.
	
	For $v_i\in V_i\setminus\{e_1^+\}$, $i=1,2$, and $v_3\in V_i\setminus\{e_1^+\}$, set $\pi_i$ to be the simple path from $e_i^-$ to $v_i$ not visiting $e_1^+$. Assume that there is an edge connecting $v_i,v_j$, $i=1,2$, $j\neq i$. If $v_j\in \pi_i$, then there is a simple path from $e_i^-$ to $v_j$ not visiting $e_1^+$, implying $v_j\in V_i$. Otherwise $v_j\notin \pi_i$, then $\pi v_j$ is a simple path from $e_i^-$ to $v_j$ not visiting $e_i^+$, which also implies $v_j\in V_i$. This is contradict to $v_j\in V_j\setminus\{e_1^+\}$ and  $V_i \cap V_j=\{e_1^+\}$.
\end{itemize}

By these three properties, we shows that
\begin{itemize}
	\item $e_1\in E_1$, $e_2\in E_2\cup E_3$;
	\item $V_1\cap(V_2\cup V_3)=\{e_1^+\}$, $E_1\cap(E_2\cup E_3)=\emptyset$;
	\item $V_1\cup(V_2\cup V_3)=V$, $E_1\cup(E_2\cup E_3)=E$.
\end{itemize}
This completes the proof.\qed

\section{Proof of Theorem \ref{thm nc complete graph}}\label{sec nc complete graph}

\noindent By Theorem \ref{thm nc adjacent edges}, we only need to show the negative correlation for disjoint edges $e_1,e_2$ in complete graph $K_n$. Consider the highest term of
\[
\mu[e_1]\mu[e_2]-\mu[e_1e_2]\mu[1].
\]
We observe that the highest term $(T[e_1]T[e_2]-T[e_1e_2]T[1])\beta^{2|V|-2}$ is vanishing. Because for any adjacent vertices $u,v$, the unit current flow $i$ from $e_1^-$ to $e_1^+$ is
\[
i(\overrightarrow{uv})=\left\{
\begin{aligned}
&2/n, &u=e_1^-,v=e_1^+,\\
&1/n, &u=e_1^-,v\neq e_1^+, \text{ or }u\neq e_1^-,v=e_1^+\\
&0, &u,v\notin \{ e_1^-,e_1^+ \}.
\end{aligned}
\right.
\] 
By Lemma \ref{lemma R decreases in r} and Proposition \ref{Prop R-eff},
\[
T[e_1]T[e_2]-T[e_1e_2]T[1]=T[e_2]T[1](R_{\rm eff}(G)-R_{\rm eff}[G/e_2])=0.
\]
Therefore, in order to verify the negative correlation, we consider the second highest term of $\mu[e_1]\mu[e_2]-\mu[e_1e_2]\mu[1]$.\\

To give a more precise estimate, we use the following property.

\begin{Prop}[Kirchhoff's Matrix-tree Theory, \cite{CJSSS2004}]
	Let $L$ be the Laplacian matrix for the graph $G$ defined by
	\[
	L_{ij}=\left\{
	\begin{aligned}
	&-c(ij), & i\neq j,\\
	&\sum_{k\neq i}c(ik), &i=j.
	\end{aligned}
	\right.
	\]
	Let $L(i)$ be the matrix given by deleting $i$-th row and column of $L$. Then for any $i$, the determinate of $L(i)$ can be expressed as
	\[
	\det L(i)=\sum_{T\ \text{\rm  is a spanning tree of }G}\left(\prod_{e\in T}c(e)\right).
	\]
\end{Prop}

By Kirchhoff's Matrix-tree Theory, we give the following lemma to show the number of spanning trees on complete graph $K_n$.

\begin{Lemma}\label{lemma number of trees}
	Set $T_n(A)$ to be the number of spanning trees on $K_n$ satisfying event $A$. Then for any disjoint edges $e_1,e_2$ in $K_n$,
	\[
	T_n[1]=n^{n-2},\ T_n[e_1]=T_n[e_2]=2n^{n-3},\ T_n[e_1e_2]=4n^{n-4}.
	\]
\end{Lemma}

\pf By Kirchhoff's Matrix-tree Theory, we compute out $T_n[1]$ directly. For $T_n[e_1]$ (this equals $T_n[e_2]$ by symmetry), there is a bijection from spanning trees on $K_n$ containing $e_1$ to spanning trees on $K_n/e_1$. Here for any positive integer $k$, each multiple edge with multiplicity $k$ in $K_n/e_1$ can be regarded as an edge with conductance $k$. Precisely, when we consider the number of spanning trees, $K_n/e_1$ can be regarded as $K_{n-1}$ with conductance
\[
c(uv)=\left\{
\begin{aligned}
&2, &u=v_0,\ \text{or }v=v_0,\\
&1, &\text{otherwise},
\end{aligned}
\right.
\]
where $v_0$ is the vertex contracted from edge $e_1$.
Based on Kirchhoff's Matrix-tree Theory, we compute out $T_n[e_1]=2n^{n-3}$.

For $T_n[e_1e_2]$, similarly, we consider $K_n/{e_1,e_2}$, which can be regarded as $K_{n-2}$ with conductance
\[
c(uv)=\left\{
\begin{aligned}
&4, &\{u,v\}=\{v_1,v_2\},\\
&2, &u\in\{v_1,v_2\},v\neq v_1,v_2,\ \text{or }v\in\{v_1,v_2\},u\neq v_1,v_2,\\
&1, &\text{otherwise},
\end{aligned}
\right.
\]
where $v_i$ is the vertex contracted from edge $e_i$ for $i=1,2$. By Kirchhoff's Matrix-tree Theory,
we compute out $T_n[e_1e_2]=4n^{n-4}$.
\qed\\

By Lemma \ref{lemma number of trees}, we verify that the second highest term of $\mu[e_1]\mu[e_2]-\mu[e_1e_2]\mu[1]$ is positive for disjoint edges $e_1,e_2$ in $K_n$.

\begin{Lemma}\label{lemma estimate 2st term for complete graph}
	For disjoint edges $e_1,e_2$ in $K_n$, if $n$ is large enough, then the coefficient of $\beta^{2|V|-3}$ in $\mu[e_1]\mu[e_2]-\mu[e_1e_2]\mu[1]$ is positive.
\end{Lemma}

By Lemma \ref{lemma estimate 2st term for complete graph}, we prove Theorem \ref{thm nc complete graph}.\\

\noindent \emph{Proof of Theorem \ref{thm nc complete graph}. }
For adjacent two edges $e_1,e_2$, by Theorem \ref{thm nc adjacent edges}, we get the negative correlation for large enough $\beta$.
For disjoint two edges $e_1,e_2$, consider $\mu[e_1]\mu[e_2]-\mu[e_1e_2]\mu[1]$. By Lemma \ref{lemma number of trees}, we deduce that the highest term should be
\[
(T_n[e_1]T_n[e_2]-T_n[e_1e_2]T_n[1])\beta^{2|V|-2}=(4n^{2n-6}-4n^{2n-6})\beta^{2|V|-2}=0.
\]
By Lemma \ref{lemma estimate 2st term for complete graph}, the second highest term should be positive for sufficiently large $n$, which implies the negative correlation for sufficiently large $\beta$.\qed\\

Here we give the proof of Lemma \ref{lemma estimate 2st term for complete graph}.\\

\noindent \emph{Proof of Lemma \ref{lemma estimate 2st term for complete graph}. }
First we consider the second highest term of $\mu[S]$ for any edge set $S$. This value should be
\[
\frac{1}{2}\sum_{V'\subset V, V'\neq V,\emptyset} F_{V'}[S]\cdot \beta^{|V|-2}.
\]
Here $F_{V'}[S]$ is the number of spanning forests containing $S$ such that there is exactly two trees $T,T'$ in the forest, where $T$ is the spanning tree of $V'$, and $T'$ is the spanning tree of $V\setminus V'$.

We then estimate the coefficient of $\beta^{2|V|-3}$ in  $\mu[e_1]\mu[e_2]-\mu[e_1e_2]\mu[1]$, which equals
\begin{align*}
&\frac{1}{2} \sum_{V'\subset V, V'\neq V,\emptyset}T_n(e_1)F_{V'}(e_2)+T_n(e_2)F_{V'}(e_1)-T_n(e_1e_2)F_{V'}(1)-T_n(1)F_{V'}(e_1e_2)\\
&=\sum_{k=1}^{[n/2]}\sum_{V'\subset V, |V'|=k}T_n(e_1)F_{V'}(e_2)+T_n(e_2)F_{V'}(e_1)-T_n(e_1e_2)F_{V'}(1)-T_n(1)F_{V'}(e_1e_2).
\end{align*}

Here we write $T_n(e_1)F_{V'}(e_2)+T_n(e_2)F_{V'}(e_1)-T_n(e_1e_2)F_{V'}(1)-T_n(1)F_{V'}(e_1e_2)$ as $a_{V'}$
For the term with $|V'|=k$ in this sum, denoted by $a_k:=\sum_{V'\subset V, |V'|=k}a_{V'}$, we compute it in the following different cases:
\begin{enumerate}
	\item the endpoints of $e_1,e_2$ are in $V'$, then the sum of $a_{V'}$ equals $-\binom{n-4}{k-4}4n^{n-4}k^{k-4}(n-k)^{n-k}$;
	\item the endpoints of $e_1,e_2$ are not in $V'$, then the sum of $a_{V'}$ equals $-\binom{n-4}{k}4n^{n-4}k^k(n-k)^{n-k-4}$;
	\item all the endpoints of $e_1$, and one endpoint of $e_2$ are in $V'$, then the sum of $a_{V'}$ equals
	
	$\binom{n-4}{k-3}8n^{n-4}k^{k-3}(n-k)^{n-k-1}$;
	\item all the endpoints of $e_2$, and one endpoint of $e_1$ are in $V'$, then the sum of $a_{V'}$ equals
	
	$\binom{n-4}{k-3}8n^{n-4}k^{k-3}(n-k)^{n-k-1}$;
	\item only one endpoint of $e_1$ is in $V'$, and the endpoints of $e_2$ are not in $V'$, then the sum of $a_{V'}$ equals $\binom{n-4}{k-1}8n^{n-4}k^{k-1}(n-k)^{n-k-3}$;
	\item only one endpoint of $e_2$ is in $V'$, and the endpoints of $e_1$ are not in $V'$, then the sum of $a_{V'}$ equals $\binom{n-4}{k-1}8n^{n-4}k^{k-1}(n-k)^{n-k-3}$;
	\item the endpoints of $e_1$ are in $V'$, and the endpoints of $e_2$ are not in $V'$, then the sum of $a_{V'}$ equals $-\binom{n-4}{k-2}4n^{n-4}k^{k-2}(n-k)^{n-k-2}$;
	\item the endpoints of $e_2$ are in $V'$, and the endpoints of $e_1$ are not in $V'$, then the sum of $a_{V'}$ equals $-\binom{n-4}{k-2}4n^{n-4}k^{k-2}(n-k)^{n-k-2}$;
	\item one endpoint of $e_1$, and one endpoint of $e_2$ are in $V'$, then the sum of $a_{V'}$ equals
	
	$-\binom{n-4}{k-2}16n^{n-4}k^{k-2}(n-k)^{n-k-2}$.
\end{enumerate}
Based on these computations, we obtain $a_k$ by summing over all values:
\[
a_k=12n^{n-3}\binom{n-4}{k-1}k^{k-4}(n-k)^{n-k-4}\frac{-k(n+6)(n-k)+2n^2}{(n-k-1)(n-k-2)}.
\]
For $k=1$, $a_1=12n^{n-3}(n-1)^{n-5}$. For $k\ge 2$, $a_k$ can be also written as
\[
-12n^{n-3}(n-1)^{n-5}\cdot\left[ \frac{(n-4)!}{(k-1)!(n-k-1)!}k^{k-4}(n-k)^{n-k-4}\frac{k(n+6)(n-k)-2n^2}{(n-1)^{n-5}} \right].
\]
For convenience, set
\[
I_k=\frac{(n-4)!}{(k-1)!(n-k-1)!}k^{k-4}(n-k)^{n-k-4}\frac{k(n+6)(n-k)-2n^2}{(n-1)^{n-5}}.
\]
We prove that for $n$ large enough, $\sum_{k=2}^{[n/2]}I_k<1$.

By Stirling's approximation, for sufficiently large $n$ and $k=pn$ for $p\in[2/n,\frac{1}{2}]$,

\begin{align*}
I_k&\le \frac{(n-4)!}{(n-k-1)!(k-1)!}k^{k-3}(n-k)^{n-k-3}(n+6)\frac{1}{(n-1)^{n-5}}\\
&\sim \frac{\sqrt{2\pi (n-4)}}{\sqrt{2\pi (n-k-1)}\sqrt{2\pi(k-1)}}\cdot\frac{(n-4)^{n-4}}{(n-k-1)^{n-k-1}(k-1)^{k-1}}\cdot e^2k^{k-3}(n-k)^{n-k-3}\frac{n+6}{(n-1)^(n-5)}\\
&=\frac{e^2}{\sqrt{2\pi}}\cdot \frac{(n-4)^{n-\frac{7}{2}}}{(n-1)^{n-\frac{7}{2}}}\cdot \frac{(n-k)^{n-k-\frac{1}{2}}}{(n-k-1)^{n-k-\frac{1}{2}}}\cdot \frac{k^{k-\frac{1}{2}}}{(k-1)^{k-\frac{1}{2}}}\cdot\frac{(n-1)^{\frac{3}{2}}(n+6)}{(n-k)^{\frac{5}{2}}}\cdot\frac{1}{k^{\frac{5}{2}}}\\
&\sim \frac{e}{\sqrt{2\pi}}\cdot\frac{(n-1)^{\frac{3}{2}}(n+6)}{(n-k)^{\frac{5}{2}}}\cdot\frac{1}{k^{\frac{5}{2}}}.
\end{align*}

Note that
\begin{align*}
\sum_{k=2}^{[n/2]}\frac{1}{[k(n-k)]^{\frac{5}{2}}}
&\le \int_{1}^{\frac{n}{2}} \frac{1}{[x(n-x)]^{\frac{5}{2}}}dx\\
&\le \int_{\frac{1}{n}}^{\frac{1}{2}} \frac{1}{[n^2p(1-p)]^{\frac{5}{2}}}dnp\\
&=\frac{16}{n^4}\cdot\left[\frac{2}{3}\cdot\frac{n-2}{2(n-1)^{\frac{1}{2}}}+\frac{1}{3}\cdot\frac{n^2(n-2)}{8(n-1)^{\frac{3}{2}}}\right].
\end{align*}

By this estimate, for $n$ sufficiently large,
\[
\sum_{k=2}^{[n/2]}I_k\le \frac{e}{\sqrt{2\pi}}\cdot\frac{2}{3}+o(1)<0.8+o(1).
\]
Since the coefficient of $\beta^{2|V|-3}$ is $12n^{n-3}(n-1)^{n-5}(1-\sum_{k=2}^{[n/2]}I_k)$, we complete our proof.\qed

\section{Proof of Theorem \ref{NG for same edges}}\label{sec nc simplify graph}

Before the proof of our main theorem, we post some lemmas.

\begin{Lemma}\label{lemma delete pivotal edges}
	If $e$ is pivotal for $G$, then negative correlation on $G\setminus \{e\}$ implies negative correlation on $G$.
\end{Lemma}

\pf Assume that the negative correlation holds on $G\setminus\{e\}$, i.e.,
\begin{equation}
\mu[e_1e_2\bar{e}]\mu[\bar{e}]\le \mu[e_1\bar{e}]\mu[e_2\bar{e}] {\rm \ for\ distinct\ }e_1,e_2\neq e.\label{1 in pivotal}
\end{equation}
Set $\mu'[S]=\mu[S \bar{e}]$. Since $e$ is pivotal for $G$, whether $e$ exists or not does not influence the existence of cycle in $G$.
Hence, $\mu[S\bar{e}]=\mu[Se]$ for $e\notin S$, by which we have
\[
\mu[S]=\left\{
\begin{aligned}
&\mu[S\bar{e}]+\beta_e\mu[Se]=(1+\beta_e)\mu'[S], &e\notin S,\\
&\mu[(S\setminus\{e\}) e]=\mu'[S\setminus\{e\}], &e\in S.
\end{aligned}
\right.
\]
Combined with (\ref{1 in pivotal}), we verify that
\[
\mu[e_1e_2]\mu[1]\le \mu[e_1]\mu[e_2] {\rm\ for\ all\ }e_1,e_2\in E.
\]
\qed

\begin{Lemma}\label{lemma every component}
	If the negative correlation holds on each component of $G$, then it holds on $G$.
\end{Lemma}

\pf We only need to prove negative correlation for $e_1,e_2$ in different components of $G$ since the existence of cycle in a component cannot be influenced by other components.

Assume that $e_i\in E_i$ for $i=1,2$, where $\{G_i=(V_i,E_i):i=1,\dots,d\}$ are all components of $G$. Set $\mu_i$ to be the measure on $G_i$ for $i=1,\dots,d$, then we deduce that ($j=1,2$)
\begin{align*}
\mu[e_1e_2]&=\mu_1[e_1]\mu_2[e_2]\prod_{i=3}^d \mu_i[1],\\
\mu[e_j]&=\mu_j[e_j]\prod_{i\neq j}\mu_i[1],\\
\mu[1]&=\prod_{i=1}^d\mu_i[1].
\end{align*}
This implies that $\mu[e_1e_2]\mu[1]\le\mu[e_1]\mu[e_2]$.  \qed

\begin{Lemma}\label{lemma delete 2-d vertex}
	The negative correlation holds on $G=(V,E)$ for any parameter $\{\beta_e\}_{e\in E}$ if and only if it holds on $\tilde{G}=(\tilde{V},\tilde{E})$ for any parameter $\{\tilde{\beta}_{\tilde{e}}\}_{\tilde{e}\in\tilde{E}}$, where $\tilde{E}$ is given by Definition \ref{def delete 2-d vertices}.
\end{Lemma}

\begin{Lemma}\label{lemma delete multiple edges}
	The negative correlation holds on a complex graph $G$ for any parameter $\{\beta_e\}_{e\in E}$ if and only if it holds on the simple graph $G':=g(G)$ for any $\{\beta_{e'}'\}_{e'\in g(E)}$, where $g$ is defined in Definition \ref{def delete multiple edges}.
\end{Lemma}

The proof of Lemmas \ref{lemma delete 2-d vertex} and \ref{lemma delete multiple edges} will be stated in the next subsection. Now we prove Theorem \ref{NG for same edges}.\\

\noindent \emph{Proof of Theorem \ref{NG for same edges}.} By Lemmas \ref{lemma delete pivotal edges} and \ref{lemma every component}, we know that the negative correlation on $G$ is equivalent to that on each component of a new graph $\hat{G}$ given by deleting all pivotal edges. Thus if $e_1$ or $e_2$ is pivotal edges, or $e_1,e_2$ belong to different component of $\hat{G}$, the negative correlation for $e_1,e_2$ holds. Without loss of generality, assume $\hat{G}$ is connected. 

By Lemmas \ref{lemma delete 2-d vertex} and \ref{lemma delete multiple edges}, we know that the negative correlation on $\hat{G}$ is equivalent to that on $g\circ f(\hat{G})$, which immediately shows the negative correlation for $e_1,e_2$ with $g\circ f(e_1)=g\circ f(e_2)$.\qed\\

Finally we give our proof of Lemmas \ref{lemma delete 2-d vertex} and \ref{lemma delete multiple edges}.\\

\noindent \emph{Proof of Lemma \ref{lemma delete 2-d vertex}.} Recall the definition of $f$ in Definition \ref{def delete 2-d vertices}. For some subgraph $F$ with edge set $E(F)$, set $\tilde{F}$ to be a graph induced by the following edge set
$\tilde{E}\setminus f(E\setminus E(F))$. We write the edge set of $\tilde{F}$ as $\tilde{E}(\tilde{F})$.
Intuitively, for any edge $e\in E$, assume $e\in f^{-1}(\tilde{e})$ for some $\tilde{e}\in \tilde{E}$, if $e\notin F$, then $\tilde{e}\notin \tilde{E}(\tilde{F})$. We show that $F\in\mathcal{F}$ is equivalent to that $\tilde{F}$ is a forest on $\tilde{G}$.

If $F\in\mathcal{F}$, set $f^{-1}(\tilde{F})$ to be a graph induced by edge set $\{e\in E: f(e)\in \tilde{E}(\tilde{F}) \}$. Note that $f^{-1}(\tilde{F})\subset F$. If there is a cycle in $\tilde{F}$, there should also exist a cycle in $f^{-1}(\tilde{F})$, contradict to $f^{-1}(\tilde{F})\subset F\in\mathcal{F}$.

If $\tilde{F}$ is a forest, assume there is a simple cycle $\pi$ in $F$. For any $e\in \pi$, we claim that
\[f^{-1}(f(e))\subset \pi.\]
By this claim we know that $f(\pi)\in \tilde{F}$. Since $f(\pi)$ is also a cycle, it is contradict to that $\tilde{F}$ is a forest.

Here we prove the claim. Otherwise, there are two adjacent edges $e_1,e_2\in f^{-1}(f(e))$ with $e_1\in\pi$ and $e_2\notin \pi$. Setting $v$ to be the vertex adjacent to $e_1,e_2$, we find that there is no simple cycle crossing $v$ in $F$ because the degree of $v$ is $1$, contradict to $e_1\in \pi$.\\

Set $\tilde{\mu}$ to be the measure of Arboreal Gas on $\tilde{G}$, where the corresponding parameter
\[
\tilde{\beta}_{\tilde{e}}=\frac{\prod_{e\in f^{-1}(\tilde{e})}\beta_e}{\prod_{e\in f^{-1}(\tilde{e})}(1+\beta_e)-\prod_{e\in f^{-1}(\tilde{e})}\beta_e}.
\]
Then for any forest $\tilde{F}'$ in $\tilde{G}$, \[\mu(\{F\in\mathcal{F}:\tilde{F}=\tilde{F}'\})=\tilde{\mu}(\tilde{F}')\cdot \prod_{\tilde{e}\in \tilde{E}}\left[\prod_{e\in f^{-1}(\tilde{e})}(1+\beta_e)-\prod_{e\in f^{-1}(\tilde{e})}\beta_e\right].\]
Here we write $\prod_{\tilde{e}\in \tilde{E}}\left[\prod_{e\in f^{-1}(\tilde{e})}(1+\beta_e)-\prod_{e\in f^{-1}(\tilde{e})}\beta_e\right]$ as $C=C(\beta_e, e\in E)$.\\

\noindent $\Leftarrow$: Consider two distinct edges $e_1,e_2\in E$.\\
\noindent \textbf{(a)}. 
If $f(e_1)=f(e_2)=\tilde{e}_0$, set
\[
\tilde{\mu}[\tilde{e}_0]=\tilde{\beta}_{\tilde{e}_0}\cdot F_1(\tilde{\beta}_{\tilde{e}},\tilde{e}\neq\tilde{e}_0),\  \tilde{\mu}[\bar{\tilde{e}}_0]= F_2(\tilde{\beta}_{\tilde{e}},\tilde{e}\neq\tilde{e}_0),
\]
where $F_2(\tilde{\beta}_{\tilde{e}},\tilde{e}\neq\tilde{e}_0)\ge F_1(\tilde{\beta}_{\tilde{e}},\tilde{e}\neq\tilde{e}_0)>0$. To verify the negative correlation on $G$, we need to prove
$
\mu[e_1]\mu[e_2]-\mu[e_1e_2]\mu[1]\ge 0.
$
This is equivalent to 
\begin{equation}
\mu[e_1\bar{e_2}]\mu[\bar{e_1}e_2]-\mu[e_1e_2]\mu[\bar{e_1}\bar{e_2}]\ge0.\label{eq nc pf in deleting 2-d vertex}
\end{equation}
The right hand side of \eqref{eq nc pf in deleting 2-d vertex}  equals
\begin{align*}
&C^2 \cdot \beta_{e_1}\prod_{e\in f^{-1}(\tilde{e}_0),e\neq e_1,e_2}(1+\beta_e)\ F_2 \cdot \beta_{e_2}\prod_{e\in f^{-1}(\tilde{e}_0),e\neq e_1,e_2}(1+\beta_e)\ F_2\\
&- C^2 \cdot \beta_{e_1}\beta_{e_2}\left[ \left( \prod_{e\in f^{-1}(\tilde{e}_0),e\neq e_1,e_2}(1+\beta_e)-\prod_{e\in f^{-1}(\tilde{e}_0),e\neq e_1,e_2}\beta_e \right)\ F_2 +\prod_{e\in f^{-1}(\tilde{e}_0),e\neq e_1,e_2}\beta_e\ F_1 \right] \\
&\hskip 2cm \cdot\prod_{e\in f^{-1}(\tilde{e}_0),e\neq e_1,e_2}(1+\beta_e)\ F_2\\
&=C^2\beta_1\beta_2\prod_{e\in f^{-1}(\tilde{e}_0),e\neq e_1,e_2}(1+\beta_e) \prod_{e\in f^{-1}(\tilde{e}_0),e\neq e_1,e_2}\beta_e\ F_2(F_2-F_1).
\end{align*}
By $F_2\ge F_1>0$, we finish the proof of \eqref{eq nc pf in deleting 2-d vertex}.\\

\noindent \textbf{(b)}. If $f(e_1)=\tilde{e}_1\neq \tilde{e}_2=f(e_2)$, set
\begin{align*}
&\tilde{\mu}[\tilde{e}_1\tilde{e}_2]=\tilde{\beta}_{\tilde{e}_1}\tilde{\beta}_{\tilde{e}_2}F_{12}(\tilde{\beta}_{\tilde{e}},\tilde{e}\neq \tilde{e}_1,\tilde{e}_2),\ &\tilde{\mu}[\tilde{e}_1\bar{\tilde{e}}_2]=\tilde{\beta}_{\tilde{e}_1}F_{1\bar{2}}(\tilde{\beta}_{\tilde{e}},\tilde{e}\neq \tilde{e}_1,\tilde{e}_2),\\
&\tilde{\mu}[\bar{\tilde{e}}_1\tilde{e}_2]=\tilde{\beta}_{\tilde{e}_2}F_{\bar{1}2}(\tilde{\beta}_{\tilde{e}},\tilde{e}\neq \tilde{e}_1,\tilde{e}_2),\ &\tilde{\mu}[\bar{\tilde{e}}_1\bar{\tilde{e}}_2]=F_{\bar{1}\bar{2}}(\tilde{\beta}_{\tilde{e}},\tilde{e}\neq \tilde{e}_1,\tilde{e}_2).
\end{align*}
By the negative correlation on $\tilde{G}$, we have
\[
F_{1\bar{2}}F_{\bar{1}2}-F_{12}F_{\bar{1}\bar{2}}\ge 0.
\]
To verify the negative correlation on $G$, we also prove \eqref{eq nc pf in deleting 2-d vertex}.
For convenience, for $i=1,2$, we rewrite $\prod_{e\in f^{-1}(\tilde{e}_i),e\neq e_i}\beta_e$, $\prod_{e\in f^{-1}(\tilde{e}_i),e\neq e_i}(1+\beta_e)$ as $C_1(e_i)$ and $C_2(e_i)$ respectively. Note the right hand side of \eqref{eq nc pf in deleting 2-d vertex} equals
\begin{align*}
&C^2\cdot \beta_{e_1}\Bigg[ \left[ C_2(e_1) - C_1(e_1) \right]C_2(e_2) F_{\bar{1}\bar{2}} + C_1(e_1)C_2(e_2) F_{1\bar{2}} \Bigg]\cdot \beta_{e_2}\Bigg[ \left[ C_2(e_2)-C_1(e_2) \right]C_2(e_1) F_{\bar{1}\bar{2}} + C_1(e_2)C_2(e_1) F_{\bar{1}2} \Bigg]\\
&-C^2\cdot\beta_{e_1}\beta_{e_2}\Bigg[ [C_2(e_1)-C_1(e_1)][C_2(e_2)-C_1(e_2)]F_{\bar{1}\bar{2}}+C_1(e_1)[C_2(e_2)-C_1(e_2)]F_{1\bar{2}}\\
&+C_1(e_2)[C_2(e_1)-C_1(e_1)]F_{\bar{1}2}+C_1(e_1)C_1(e_2)F_{12} \Bigg]\cdot C_2(e_1)C_2(e_2)F_{\bar{1}\bar{2}}\\
&=C^2\prod_{i=1,2} \beta_{e_i} C_1(e_i)C_2(e_i) (F_{1\bar{2}}F_{\bar{1}2}-F_{12}F_{\bar{1}\bar{2}}),
\end{align*}
which is non-negative because $F_{1\bar{2}}F_{\bar{1}2}-F_{12}F_{\bar{1}\bar{2}}\ge 0$. This completes the proof of \eqref{eq nc pf in deleting 2-d vertex}.\\

\noindent $\Rightarrow$: For distinct $\tilde{e}_1,\tilde{e}_2\in\tilde{E}$, by \textbf{(b)} in ``$\Leftarrow$", the negative correlation on $G$ is equivalent to
\[
F_{1\bar{2}}F_{\bar{1}2}-F_{12}F_{\bar{1}\bar{2}}\ge 0.
\]
This shows immediately that
\[
\tilde{\mu}[\tilde{e}_1\bar{\tilde{e}}_2]\tilde{\mu}[\bar{\tilde{e}}_1\tilde{e}_2]-\tilde{\mu}[\tilde{e}_1\tilde{e}_2]\tilde{\mu}[\bar{\tilde{e}}_1\bar{\tilde{e}}_2]\ge0,
\]
which completes the proof of negative correlation on $\tilde{G}$.\qed\\

\noindent \emph{Proof of Lemma \ref{lemma delete multiple edges}.} Recall that $g$ maps edge sets to edge sets by deleting multiple edges. For some subgraph $F=(V(F),E(F))$ in $G$, denote by $F'$ the graph induced by edge set $E'(F'):=g(E(F))$. If $F$ is a forest, there should not be any multiple edges in $F$. This implies that $F$ has the same structure as that of $F'$. On the other hand, if $F'$ is a forest, set $g^{-1}(F')$ to be a graph induced by edge set $\{e\in E: g(e)\in E'(F') \}$. Then all the subgraph of $g^{-1}(F')$ without multiple edges is a forest.

Set $\mu'$ to be the measure of Arboreal Gas on $G'$ with the parameter
\[
\beta_{e'}'=\sum_{e\in g^{-1}(e')}\beta_e.
\]
Then for any forest $F''$ in $G'$,
\[
\mu(\{F\in\mathcal{F}: F'=F'' \})=\mu'(F'').
\]

\noindent $\Leftarrow$: Consider two distinct edges $e_1,e_2\in E$.\\
\noindent \textbf{(a)}. If $e_1,e_2\in g^{-1}(e_0')$ for some $e_0'\in E'$, set
\[
\mu'[e_0']=\beta_{e_0'}'F_1(\beta_{e'}',e'\neq e_0'),\ \mu'[\bar{e_0'}]=F_2(\beta_{e'}',e'\neq e_0'),
\]
where $F_2\ge F_1>0$. The negative correlation on $G$ is equivalent to \begin{equation}
\mu[e_1\bar{e_2}]\mu[\bar{e_1}e_2]-\mu[e_1e_2]\mu[\bar{e_1}\bar{e_2}]\ge0.\label{eq nc pf in deleting multiple edges}
\end{equation}
Note that $\mu[e_1e_2]=0$. The right hand side of \eqref{eq nc pf in deleting multiple edges} is equal to $\mu[e_1\bar{e_2}]\mu[\bar{e_1}e_2]\ge0$.\\

\noindent \textbf{(b)}. If $g(e_1)=e_1'\neq e_2'=g(e_2)$, set
\begin{align*}
&\mu'(e_1'e_2')=\beta_{e_1'}'\beta_{e_2'}'F_{12}(\beta_{e'}',e'\neq e_1',e_2'),\ &\mu'(e_1'\bar{e_2'})=\beta_{e_1'}'F_{1\bar{2}}(\beta_{e'}',e'\neq e_1',e_2'),\\
&\mu'(\bar{e_1'}e_2')=\beta_{e_2'}'F_{\bar{1}2}(\beta_{e'}',e'\neq e_1',e_2'),\ &\mu'(\bar{e_1'}\bar{e_2'})=F_{\bar{1}\bar{2}}(\beta_{e'}',e'\neq e_1',e_2').
\end{align*}

By the negative correlation on $G'$, we obtain
\[
F_{1\bar{2}}F_{\bar{1}2}-F_{12}F_{\bar{1}\bar{2}}\ge 0.
\]
Now we verify \eqref{eq nc pf in deleting multiple edges} to prove the negative correlation on $G$. Rewrite $\prod_{e\in g^{-1}(e_i),e\neq e_i}\beta_e$ as $C(e_i)$ for $i=1,2$. Then the right hand side of \eqref{eq nc pf in deleting multiple edges} should be
\begin{align*}
&\beta_{e_1}(C(e_2)F_{12}+F_{1\bar{2}})\cdot \beta_{e_2}(C(e_1)F_{12}+F_{\bar{1}2})\\
&-\beta_{e_1}\beta_{e_2}F_{12}\cdot(C(e_1)C(e_2)F_{12}+C(e_1)F_{1\bar{2}}+C(e_2)F_{\bar{1}2}+F_{\bar{1}\bar{2}})\\
&=\beta_{e_1}\beta_{e_2}(F_{1\bar{2}}F_{\bar{1}2}-F_{12}F_{\bar{1}\bar{2}}).
\end{align*}
By $F_{1\bar{2}}F_{\bar{1}2}-F_{12}F_{\bar{1}\bar{2}}\ge 0$, we deduce \eqref{eq nc pf in deleting multiple edges}.\\

\noindent $\Rightarrow$: For distinct $e_1',e_2'\in E'$, by \textbf{(b)} in ``$\Leftarrow$", the negative correlation on $G$ is equivalent to \[F_{1\bar{2}}F_{\bar{1}2}-F_{12}F_{\bar{1}\bar{2}}\ge 0,\]
which implies that
\[
\mu'[e_1'\bar{e_2'}]\mu'[\bar{e_1'}e_2']-\mu'[e_1'e_2']\mu[\bar{e_1'}\bar{e_2'}]\ge0.
\]
This shows the negative correlation on $G'$. \qed

\bibliographystyle{abbrvnat}
	\bibliography{ref.bib}

\end{document}